\documentclass[10pt]{amsart}

\usepackage{amsmath,amsthm,amsfonts,amscd,amssymb}
\usepackage{hyperref,wasysym}
\usepackage{verbatim}
\usepackage[pdftex]{graphicx}

\newtheorem{thm}{Theorem}[section]
\newtheorem{cor}[thm]{Corollary}
\newtheorem{lem}[thm]{Lemma}

\newtheorem{quest}{Question}

\theoremstyle{definition}

\theoremstyle{remark}

\begin{document}

\title{Deep hole lattices and isogenies of elliptic curves}

\author{Lenny Fukshansky}\thanks{Fukshansky was partially supported by the Simons Foundation grant \#519058}
\author{Pavel Guerzhoy}
\author{Tanis Nielsen}

\address{Department of Mathematics, Claremont McKenna College, 850 Columbia Avenue, Claremont, CA 91711}
\email{lenny@cmc.edu}
\address{Department of Mathematics, University of Hawaii, 2565 McCarthy Mall, Honolulu, HI, 96822-2273}
\email{pavel@math.hawaii.edu}
\address{Department of Mathematics, Harvey Mudd College, 320 E. Foothill Blvd., Claremont CA 91711}
\email{tnielsen@g.hmc.edu}


\begin{abstract} Given a lattice $L$ in the plane, we define the affiliated deep hole lattice $H(L)$ to be spanned by a shortest vector of $L$ and a deep hole of $L$ contained in the triangle with sides corresponding to the shortest basis vectors. We study the geometric and arithmetic properties of deep hole lattices. In particular we investigate conditions on $L$ under which $H(L)$ is well-rounded and prove that $H(L)$ is defined over the same field as $L$. For the period lattice corresponding to an isomorphism class of elliptic curves, we produce a finite sequence of deep hole lattices ending with a well-rounded lattice which corresponds to a point on the boundary arc of the fundamental strip under the action of $\operatorname{SL}_2(\mathbb{Z})$ on the upper halfplane. In the case of CM elliptic curves, we prove that all elliptic curves generated by this sequence are isogenous to each other and produce bounds on the degree of isogeny. Finally, we produce a counting estimate for the planar lattices with a prescribed deep hole lattice.
\end{abstract}

\maketitle

\def\A{{\mathcal A}}
\def\B{{\mathcal B}}
\def\C{{\mathcal C}}
\def\D{{\mathcal D}}
\def\F{{\mathcal F}}
\def\x{{\mathcal H}}
\def\I{{\mathcal I}}
\def\J{{\mathcal J}}
\def\K{{\mathcal K}}
\def\L{{\mathcal L}}
\def\M{{\mathcal M}}
\def\N{{\mathcal N}}
\def\O{{\mathcal O}}
\def\R{{\mathcal R}}
\def\s{{\mathcal S}}
\def\V{{\mathcal V}}
\def\W{{\mathcal W}}
\def\X{{\mathcal X}}
\def\Y{{\mathcal Y}}
\def\H{{\mathcal H}}
\def\Z{{\mathcal Z}}
\def\OO{{\mathcal O}}
\def\BB{{\mathbb B}}
\def\cee{{\mathbb C}}
\def\pee{{\mathbb P}}
\def\que{{\mathbb Q}}
\def\real{{\mathbb R}}
\def\zed{{\mathbb Z}}
\def\NN{{\mathbb N}}
\def\hyp{{\mathbb H}}
\def\aa{{\mathfrak a}}
\def\HH{{\mathfrak H}}
\def\qbar{{\overline{\mathbb Q}}}
\def\eps{{\varepsilon}}
\def\ahat{{\hat \alpha}}
\def\bhat{{\hat \beta}}
\def\gt{{\tilde \gamma}}
\def\h{{\tfrac12}}
\def\be{{\boldsymbol e}}
\def\bei{{\boldsymbol e_i}}
\def\bff{{\boldsymbol f}}
\def\ba{{\boldsymbol a}}
\def\bb{{\boldsymbol b}}
\def\bc{{\boldsymbol c}}
\def\bm{{\boldsymbol m}}
\def\bk{{\boldsymbol k}}
\def\bi{{\boldsymbol i}}
\def\bl{{\boldsymbol l}}
\def\bq{{\boldsymbol q}}
\def\bu{{\boldsymbol u}}
\def\bt{{\boldsymbol t}}
\def\bs{{\boldsymbol s}}
\def\bv{{\boldsymbol v}}
\def\bw{{\boldsymbol w}}
\def\bx{{\boldsymbol x}}
\def\bX{{\boldsymbol X}}
\def\bz{{\boldsymbol z}}
\def\bwy{{\boldsymbol y}}
\def\bY{{\boldsymbol Y}}
\def\bL{{\boldsymbol L}}
\def\baa{{\boldsymbol\alpha}}
\def\bbb{{\boldsymbol\beta}}
\def\bet{{\boldsymbol\eta}}
\def\btau{{\boldsymbol\tau}}
\def\bgamma{{\boldsymbol\gamma}}
\def\bxi{{\boldsymbol\xi}}
\def\bo{{\boldsymbol 0}}
\def\bol{{\boldkey 1}_L}
\def\ep{\varepsilon}
\def\p{\boldsymbol\varphi}
\def\q{\boldsymbol\psi}
\def\rank{\operatorname{rank}}
\def\aut{\operatorname{Aut}}
\def\lcm{\operatorname{lcm}}
\def\sgn{\operatorname{sgn}}
\def\spn{\operatorname{span}}
\def\md{\operatorname{mod}}
\def\Norm{\operatorname{Norm}}
\def\dim{\operatorname{dim}}
\def\det{\operatorname{det}}
\def\Vol{\operatorname{Vol}}
\def\rk{\operatorname{rk}}
\def\Gal{\operatorname{Gal}}
\def\WR{\operatorname{WR}}
\def\WO{\operatorname{WO}}
\def\GL{\operatorname{GL}}
\def\SL{\operatorname{SL}}
\def\pr{\operatorname{pr}}
\def\Tr{\operatorname{Tr}}

\section{Introduction}
\label{intro}

Let $L \subset \real^2$ be a lattice with successive minima $\lambda_1 \leq \lambda_2$ and the corresponding minimal basis vectors $\bx_1,\bx_2$. It is well known that, choosing $\pm \bx_1,\pm \bx_2$ if necessary, we can ensure that the angle $\theta$ between these vectors is in the interval $[\pi/3,\pi/2]$: this angle is an invariant of the lattice, we call it {\it angle} of $L$. The lattice $L$ is called {\it well-rounded} (WR) if $\lambda_1 = \lambda_2$; it is called {\it semi-stable} if $\lambda_1 \geq \det(L)^{1/2}$. For planar lattices, WR property implies semi-stability (this is not true in higher dimensions). Two lattices $L_1$ and $L_2$ are said to be {\it similar}, denoted $L_1 \sim L_2$, if there exists $\alpha \in \real_{>0}$ and $U \in \O_2(\real)$ such that $L_2 = \alpha U L_1$. This is an equivalence relation which preserves WR and semi-stability properties, hence we will speak of WR and semi-stable similarity classes. WR and semi-stable lattices play an important role in lattice theory and are objects of substantial study; see \cite{casselman}, \cite{lf:pg:fl}, \cite{lf:dk} and references within for some details.

\begin{figure}[t]
\centering
\includegraphics[scale=0.4]{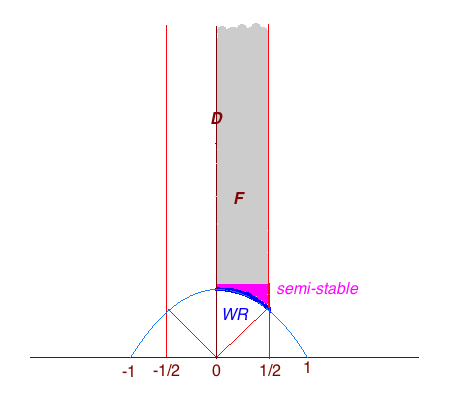}
\caption{Space of lattices in $\real^2$ with WR and semi-stable subregions marked by colors.}\label{fig:domain}
\end{figure}

Similarity classes of planar lattices can be parameterized as follows. Let $\hyp = \{ \tau = a+bi : b \geq 0 \} \subset \cee$ be the upper half-plane, and let
$$\D := \{ \tau = a + bi \in \hyp : -1/2 < a \leq 1/2, |\tau| \geq 1 \}.$$
Let
$$\F := \{ \tau = a + bi \in \hyp : 0 \leq a \leq 1/2, |\tau| \geq 1 \},$$
so, loosely speaking, $\F$ is ``half" of $\D$. Every point $\tau = a + bi \in \F$ can be identified with a lattice
\begin{equation}
\label{Ltau}
\Lambda_{\tau} := \begin{pmatrix} 1 & a \\ 0 & b \end{pmatrix} \zed^2
\end{equation}
in $\real^2$. Then the lattice $L$ with successive minima $\lambda_1 \leq \lambda_2$ and the corresponding minimal basis vectors $\bx_1, \bx_2$ is similar to $L'= \frac{1}{\lambda_1}L$. Rotating if necessary, we can ensure that the image of $\bx_1$ under this similarity coincides with $\be_1 := \begin{pmatrix} 1 \\ 0 \end{pmatrix}$. Then the image $\bx'_2$ of $\bx_2$ must have its first coordinate between $0$ and $1/2$, since otherwise one of $\bx'_2 \pm \be_1 \in L$ would be a shorter vector than $\bx'_2$ and still linearly independent with $\be_1$. Furthermore, reflecting over $\be_1$, if necessary, we can assume that $\bx'_2$ has a positive second coordinate. In other words, every planar lattice $L$ is similar to a  lattice of the form $\Lambda_{\tau}$ for some $\tau \in \F$, and it is a well-known fact that it is similar to {\it precisely one} such lattice; we will say that this similarity class is {\it represented} by $\tau$. Hence, $\F$ can be thought of as the space of similarity classes of lattices in~$\real^2$ (see Figure~\ref{fig:domain}). WR similarity classes correspond to the circular arc $\{ \tau \in \F : |\tau|=1 \}$ and semi-stable similarity classes correspond to the set $\{ \tau = a+bi \in \F : b \leq 1 \}$.

On the other hand, the full domain $\D$ can be viewed as the space of isomorphism classes of elliptic curves: a point $\tau$ corresponds to the isomorphism class of the elliptic curve given by the complex torus $\cee/\Lambda_{\tau}$, where we are identifying $\cee$ with $\real^2$ and thinking of $\Lambda_{\tau}$ as $\spn_{\zed} \{1,\tau\} \subset \cee$. This being said, while the lattices $\Lambda_{\tau}$ and $\Lambda_{-\bar{\tau}}$ are similar the corresponding elliptic curves are not isomorphic: instead, the two elliptic curves have conjugate $j$-invariants, since $j(-\bar{\tau}) = \overline{j(\tau)}$  (here $j$ is Klein's modular $j$-function; see~\cite{silverman} for details).

A {\it deep hole} of a lattice $L$ is a point in~$\real^2$ which is farthest away from the lattice. The distance from the origin to the nearest deep hole is the {\it covering radius} $\mu$ of $L$. There is a unique deep hole $\bz$ of $L$ contained in the triangle $T$ with vertices $\bo$ and the endpoints of $\bx_1,\bx_2$: we call it the {\it fundamental deep hole} of $L$. Define the {\it deep hole lattice} of $L$ to be
$$H(L) := \spn_{\zed} \{ \bx_1, \bz \}.$$
The main goal of this paper is to investigate the properties of this new family of deep hole lattices and connect them to the properties of the associated elliptic curves. We start with some observations about the geometry and arithmetic of $H(L)$, proving the following theorem in Section~\ref{dhl}.

\begin{thm} \label{main1} Let $L$ be a lattice in the plane with the angle $\theta \in [\pi/3,\pi/2]$ and successive minima $\lambda_1$ and $\lambda_2 = \alpha \lambda_1$ for some $\alpha \geq 1$. Let $H(L)$ be the deep hole lattice of $L$. The following statements hold:
\begin{enumerate}

\item If $\alpha \leq 2 \sin (\theta + \pi/6)$, then $H(L)$ is WR.

\item If $L$ is semi-stable, then $H(L)$ is WR.

\item If $L$ is WR, then $H(L) \sim L$.

\item  If $L \subset K^2$ for some subfield $K$ of $\real$, then $H(L) \subset K^2$.

\end{enumerate}

\end{thm}

Next we turn our attention specifically to the lattices of the form $\Lambda_{\tau}$ for $\tau \in \F$ parameterizing all the similarity classes in the plane. Given a subfield $K$ of $\real$, we say that a similarity class represented by $\tau$ {\it lies over} $K$ if $\tau = a+bi$ with real numbers $a,b \in K$. This is equivalent to saying that some lattice in this similarity class is contained in $K(i) \subseteq \cee$, which is identified with $K^2 \subseteq \real^2$. 


A similarity class represented by $\tau$ lying over $K$ contains infinitely many lattices $L \subset K^2$. It is not difficult to observe that the similarity class of the deep hole lattice lies over the same field as the similarity class of the original lattice (Lemma~\ref{tau_dh}). In fact, the deep hole construction allows us to describe an interesting family of similarity classes defined over the same field.

\begin{thm} \label{dh_seq} Let $\tau_0 = a_0+b_0 i \in \F$ with $a_0,b_0 \in K$ for some subfield $K \subseteq \real$. There exists a finite sequence of numbers $\tau_1,\dots,\tau_n$ given by $\tau_k = a_k+b_ki$ for all $1 \leq k \leq n$, so that
\begin{equation}
\label{dh_seq1}
a_k = \frac{1}{2},\ b_k = \frac{a_{k-1}^2+b_{k-1}^2-a_{k-1}}{2b_{k-1}} \in K\ \forall\ 1 \leq k \leq n,
\end{equation}
with $\Lambda_{\tau_k} = H(\Lambda_{\tau_{k-1}})$ and $\Lambda_{\tau_n}$ well-rounded, hence $H(\Lambda_{\tau_n}) \sim \Lambda_{\tau_n}$. Furthermore, $\tau_1,\dots,\tau_{n-1} \in \F$, $|\tau_n| \leq 1$ and $n \leq \log_2 \left( \frac{2b_0}{\sqrt{3}} \right)$.
\end{thm}

We prove Theorem~\ref{dh_seq} in Section~\ref{iso}. For a given $\tau_0 = a_0+b_0i \in \F$, we will refer to the sequence of complex numbers
$$\tau_k = a_k+b_ki,\ 1 \leq k \leq n,$$
defined as in~\eqref{dh_seq1}, as the {\it deep hole sequence} for $\tau_0$ and call the corresponding sequence of lattices $\Lambda_{\tau_k}$ as defined in Corollary~\ref{dh_seq} the {\it deep hole lattice sequence} for $\Lambda_{\tau_0}$.

We focus on a particular situation in which this deep hole lattice sequence is especially interesting. A lattice $L = A\zed^2 \subset \real^2$ is called {\it arithmetic} if the corresponding Gram matrix $A^{\top} A$ is a scalar multiple of an integer matrix, i.e. $L$ is similar to an integral lattice. A lattice $\Lambda_{\tau}$ is arithmetic if and only if
$$\tau = \frac{p}{q} + i \sqrt{\frac{s}{t}} \in \F$$
for some $p,q,s,t \in \zed$, i.e. $\tau$ is a quadratic irrationality. This condition is equivalent to the elliptic curve $E_{\tau}$ with period lattice $\Lambda_{\tau}$ having complex multiplication by the imaginary quadratic field~$\que(\tau)$, i.e. the endomorphism ring of $E_{\tau}$ is an order in~$\que(\tau)$ properly containing $\zed$. Thus, if $\tau_0$ in Theorem~\ref{dh_seq} is a quadratic irrationality, then $\Lambda_{\tau_0}$ and the entire sequence of deep hole lattices $\Lambda_{\tau_1},\dots,\Lambda_{\tau_n}$ are arithmetic with the corresponding CM elliptic curves $E_{\tau_k}$ having their endomorphism rings in the same imaginary quadratic field~$\que(\tau_0)$. Furthermore, each curve $E_{\tau_k}$ for $1 \leq k \leq n-1$ has a real negative $j$-invariant, whereas the curve $E_{\tau_n}$ corresponds to a well-rounded lattice $\Lambda_{\tau_n}$, and hence has real $j$-invariant in the interval $[0,1]$ (see~\cite{lf:pg:fl} for details). We review the necessary notation and prove the following result in Section~\ref{iso}.

\begin{thm} \label{main_iso} Let $\tau_0 = a_0 + b_0 i \in \F$ be a quadratic irrationality and 
$$\{ \tau_k = a_k+b_k i \}_{k=1}^n$$
its corresponding deep hole sequence. For each $0 \leq k \leq n$, let $E_{\tau_k}$ be the corresponding CM elliptic curve with the arithmetic period lattice $\Lambda_{\tau_k}$. Then all of these elliptic curves are isogenous. Furthermore, for any $0 \leq k \leq n-1$, there exists an isogeny between $E_{\tau_k}$ and $E_{\tau_{k+1}}$ with degree
$$\delta_k \leq \frac{12\sqrt{3}\ b_{k+1}\ d_k^4\ (a_k^2+b_k^2)^2}{b_k},$$
where $d_k = \min \{ d \in \zed_{>0} : da_k, d^2b^2_k \in \zed\}$.
\end{thm}

Finally, we consider a certain inverse problem. From here on, let $K$ be an arbitrary number field of degree $n$ and suppose that the similarity class represented by $\tau_0 = \frac{1}{2} + it \in \F$ lies over $K$. Consider the set 
\begin{equation}
\label{SK}
S_{K,\tau_0} = \left\{ \tau \in \F : \tau \text{ is defined over $K$ and } H(\Lambda_{\tau}) = \Lambda_{\tau_0} \right\}.
\end{equation}
i.e., the set of similarity classes defined over $K$ whose deep hole lattice is $\Lambda_{\tau_0}$. While this is an infinite set, we can count these similarity classes bounding the so-called primitive height $\H^p$ of $\tau$ which is defined in Section~\ref{count}.

\begin{thm} \label{ht_count} For a real number $T \geq 1$, define
$$S_{K,\tau_0}(T) = \left\{ \tau \in S_{K,\tau_0} : \H^p(\tau) \leq T \right\}.$$
Then, as $T \to \infty$,
$$|S_{K,\tau_0}(T)| \leq \left( \frac{4^{r_1} \pi^{2r_2}}{8 \zeta(2n) \left( 2t+\sqrt{4t^2+1} \right) |\Delta_K|} \right) T^{2n} +  O(T^{2n-1}),$$
where $\zeta$ stands for the Riemann zeta-function and $n = [K:\que]$.
\end{thm} 

\noindent
We recall necessary notation and prove Theorem~\ref{ht_count} in Section~\ref{count}. We are now ready to proceed.
\bigskip

\section{Geometry of deep hole lattices}
\label{dhl}

The main goal of this section is to prove Theorem~\ref{main1}. We do it via a lemma and three corollaries, each proving one of the four parts of the theorem. As in the statement of the theorem, let us write $\alpha$ for the quotient $\lambda_2/\lambda_1 \geq 1$.

\begin{lem} \label{dh_WR} If $\alpha \leq 2 \sin (\theta + \pi/6)$, then $H(L)$ is WR.
\end{lem}

\proof
Let $\nu$ be the angle between the vectors $\bz,\bx_1-\bz$. Since $\bz$ is the center of the circle circumscribed about the triangle $T$, we have $\|\bz\| = \|\bx_1-\bz\|$ as both of these line segments are radii of this circle. This common value is $\mu$, the covering radius of the lattice $L$, which is given by the formula
$$\mu = \frac{\sqrt{\lambda_1^2 + \lambda_2^2 - 2\lambda_1\lambda_2\cos \theta}}{2 \sin \theta},$$
see Lemma~5.2 of~\cite{mf_lf}. 

Then the triangle with sides 
$$\bx_1,\bz,\bx_1-\bz$$
is isosceles, and so the bisector of the angle $\nu$ is perpendicular to $\bx_1$ and is the median of that side. Therefore we have
\begin{equation}
\label{sin-nu}
\sin(\nu/2) = \frac{\lambda_1/2}{\mu} = \frac{\lambda_1 \sin \theta}{\sqrt{\lambda_1^2+\lambda_2^2-2\lambda_1 \lambda_2 \cos \theta}} = \frac{\sin \theta}{\sqrt{1+\alpha^2-2 \alpha \cos \theta}}.
\end{equation}
If $\nu \in [\pi/3,\pi/2]$, then the lattice $H(L)$ is well-rounded: since $\|\bz\| = \|\bx_1-\bz\|$,
$$\| a \bz + b (\bx_1-\bz) \|^2 = (a^2+b^2) \|\bz\|^2 + 2ab \|\bz\|^2 \cos \nu \geq \|\bz\|^2,$$
for any $a,b \in \zed$, and hence $\bz,\bx_1-\bz$ is a minimal basis for $H(L)$. Thus the sufficient condition is $\sin (\nu/2) \geq \sin (\pi/6) = 1/2$, and so we have the condition:
$$4 \sin^2 \theta \geq 1 + \alpha^2 - 2\alpha \cos \theta.$$
Using the identity $\sin^2 \theta = 1-\cos^2 \theta$ and rearranging, we obtain
$$3-3 \cos^2 \theta = 3 \sin^2 \theta \geq \alpha^2 - 2\alpha \cos \theta + \cos^2 \theta = (\alpha-\cos \theta)^2,$$
which leads to
$$\alpha \leq \sqrt{3} \sin \theta + \cos \theta = 2 \sin (\theta + \pi/6).$$
\endproof

Now suppose that $L \sim \Lambda_{\tau}$ for some 
\begin{equation}
\label{tau}
\tau = \alpha e^{i\theta} = \alpha \cos \theta + \alpha i \sin \theta \in \F,
\end{equation}
where $\alpha = |\tau|$ and $\theta$ is the angle of $L$. 

\begin{cor} \label{dh_stable} If $L$ is semi-stable, then $H(L)$ is WR.
\end{cor}

\proof
If $L$ is semi-stable, then
$$\alpha = |\tau| \leq \sqrt{5}/2 < \sqrt{3} = \min \left\{ 2 \sin (\theta + \pi/6) : \theta \in [\pi/3,\pi/2] \right\},$$
and so the condition of Lemma~\ref{dh_WR} is satisfied.
\endproof

\begin{cor} \label{dh_WR_sim} If $L$ is WR, then $H(L) \sim L$.
\end{cor}

\proof
Notice that two WR lattices are similar if and only if they have the same angle. Using notation and setup of Lemma~\ref{dh_WR} and applying~\eqref{sin-nu} with $\alpha=1$ (since $L$ is WR), we obtain:
$$\sin(\nu/2) =  \frac{\sin \theta}{\sqrt{2-2 \cos \theta}},\ \cos(\nu/2) = \frac{\sqrt{2 - 2 \cos \theta - \sin^2 \theta}}{\sqrt{2-2 \cos \theta}}.$$
Therefore
\begin{eqnarray*}
\cos \nu & = & \cos^2 (\nu/2) - \sin^2(\nu/2) = \frac{2 - 2 \cos \theta - 2\sin^2 \theta}{2-2 \cos \theta} \\
& = &  \frac{- \cos \theta (1 - \cos \theta)}{1 - \cos \theta} = -\cos \theta = \cos (\pi-\theta).
\end{eqnarray*}
Hence the angle between $\bz$ and $-(\bz-\bx_1) = \bx_1-\bz$ is $\theta \in [\pi/3,\pi/2]$. Then $\bz, \bx_1-\bz$ is a minimal basis for $H(L)$, so the angle of $H(L)$ is $\theta$.
\endproof

\begin{cor} \label{same_field} If $L \subset K^2$ for some subfield $K$ of $\real$, then $H(L) \subset K^2$, and hence the similarity class of $H(L)$ lies over $K$.
\end{cor}

\proof
Since $L \subset K^2$, we have $\lambda_j^2 = \|\bx_j\|^2 \in K$ for $j=1,2$, and so $\alpha^2 \in K$. Let $A = (\bx_1\ \bx_2) \in \GL_2(K)$ be the corresponding basis matrix. By Lemma~4.1 of~\cite{mf_lf}, the fundamental deep hole of $L$ is
\begin{equation}
\label{z_dh}
\bz = \frac{1}{2} (A^{\top})^{-1} \begin{pmatrix} \lambda_1^2 \\ \lambda_2^2 \end{pmatrix},
\end{equation}
and hence $H(L) = \spn_{\zed} \left\{ \bx_1,\bz \right\} \subset K^2$. As we detailed in Section~\ref{intro}, this implies that the similarity class of $H(L)$ lies over $K$.
\endproof
\bigskip

\section{Deep hole lattices in the fundamental strip}
\label{iso}

We start with a simple observation about the similarity class of a deep hole lattice.

\begin{lem} \label{tau_dh} Let $\tau = a+bi \in \F$, $\Lambda_{\tau}$ the corresponding lattice as in~\eqref{Ltau} and suppose that $\Lambda_{\tau}$ lies over a field $K \subseteq \cee$. Then $\gamma \in \F$ representing the similarity class of $H(\Lambda_{\tau})$ also lies over $K$.
\end{lem}

\proof
Let us write
$$\be_1 = \begin{pmatrix} 1 \\ 0 \end{pmatrix},\ \btau = \begin{pmatrix} a \\ b \end{pmatrix}$$
for the minimal basis of $\Lambda_{\tau}$ and $\theta \in [\pi/3,\pi/2]$ for the angle between these two vectors, hence
$$\cos \theta = \frac{a}{\sqrt{a^2+b^2}},\ \sin \theta = \frac{b}{\sqrt{a^2+b^2}}.$$
Then~\eqref{z_dh} guarantees that the fundamental deep hole of $\Lambda_{\tau}$ is
\begin{equation}
\label{tau1}
\btau_1 = \begin{pmatrix} 1/2 \\ (a^2+b^2-a)/2b \end{pmatrix},
\end{equation}
which certainly lies over $K$ and $H(\Lambda_{\tau}) = \Lambda_{\tau_1}$ for $\tau_1 = \frac{1}{2} + \frac{a^2+b^2-a}{2b} i$. If $\tau_1 \in \F$, then $\gamma = \tau_1$ and we are done.

Suppose $\tau_1 \not\in \F$. This means that $|\tau_1| < 1$, then the vectors $\btau_1, \be_1 - \btau_1 \in \Lambda_{\tau_1}$ have equal length and the angle between them is in the interval $(\pi/3,2\pi/3]$: the largest value $2\pi/3$ is attained when $\tau = \frac{1}{2} + \frac{\sqrt{3}}{2}i$. On the other hand, $|\tau_1| < 1$ implies that $(a^2+b^2-a)/2b < \sqrt{3}/2$, hence
$$\frac{\|\btau\|}{\|\be_1\|} = \sqrt{a^2+b^2} < \frac{\sqrt{3}b + a}{\sqrt{a^2+b^2}} = \sqrt{3} \sin \theta + \cos \theta,$$
which means that $\Lambda_{\tau}$ satisfies the condition of Lemma~\ref{dh_WR}. Therefore $H(\Lambda_{\tau})$ is WR and its similarity class $\gamma$ lies over $K$ by Corollary~\ref{same_field}. This completes the proof.
\endproof

We can use this lemma to give a quick proof of Theorem~\ref{dh_seq}.

\proof[Proof of Theorem~\ref{dh_seq}]
Existence of the sequence $\tau_k$ with the claimed properties follows by iterative application of~\eqref{tau1}. Further, $a_k \leq 1/2$ for each $k$, so
$$b_k = \frac{b_{k-1}}{2} - \frac{a_{k-1} (1 - a_{k-1})}{2b_{k-1}} \leq \frac{b_{k-1}}{2} - \frac{1}{8b_{k-1}} < \frac{b_{k-1}}{2},$$
hence $b_k < \frac{b_0}{2^k}$. If $n \geq \log_2 \left( \frac{2b_0}{\sqrt{3}} \right)$, then $b_n < \sqrt{3}/2$, which means that $|\tau_n| < 1$ (i.e., $\tau_n \not\in \F$) and hence $\Lambda_{\tau_n}$ is well-rounded. Then Corollary~\ref{dh_WR_sim} implies that $H(\Lambda_{\tau_n})$ is similar to $\Lambda_{\tau_n}$ and hence the sequence of deep holes has stabilized. 
\endproof

Next we turn to the proof of Theorem~\ref{main_iso}. First, let us recall the definition of an isogeny. Given lattices $\Lambda, \Lambda' \subset \cee$ a nonzero morphism $E \rightarrow E'$ between the corresponding elliptic curves $E=\cee/\Lambda$ and $E'=\cee/\Lambda'$ which takes $0$ to $0$ is called an {\it isogeny}. An isogeny is always surjective and has a finite kernel. For instance, if $\beta \in \cee^*$ is such that $\beta \Lambda \subseteq \Lambda'$, then multiplication by $\beta$ function $z \mapsto \beta z$ maps $\cee \rightarrow \cee$ modulo the lattices, and hence maps $E \rightarrow E'$. In fact, every isogeny is of this form. The {\it degree} of this isogeny is then the degree of the morphism, which is equal to the index of the sublattice $\beta \Lambda$ in $\Lambda'$, i.e.
$$\deg(\beta) = [\Lambda':\beta \Lambda].$$
This is precisely the size of its kernel. If an isogeny $E \to E'$ exists, then there also exists the dual isogeny $E' \to E$ of the same degree such that their composition is simply the multiplication-by-degree map, and hence the curves are called isogenous: this is an equivalence relation. There may exist multiple isogenies between two elliptic curves, but since degree of an isogeny is a positive integer, we can ask for an isogeny of minimal degree. See \cite{silverman} for more details on elliptic curves and their isogenies.
\bigskip

\proof[Proof of Theorem~\ref{main_iso}]
To prove the first statement, pick any $0 \leq k \leq n-1$, then $\btau_{k+1} = \begin{pmatrix} a_{k+1} \\ b_{k+1} \end{pmatrix}$ is the fundamental deep hole of $\Lambda_{\tau_k}$. Notice that the successive minima of $\Lambda_{\tau_k}$ are $\lambda_{k1} = 1,\ \lambda_{k2} = |\tau_k|$ with the corresponding minimal basis vectors $\be_1$, $\btau_k$, and hence
$$\be_1 \cdot \btau_k = a_k \lambda_{k1}^2 = \frac{a_k}{a_k^2+b_k^2} \lambda_{k2}^2,$$
where $a_k, \frac{a_k}{a_k^2+b_k^2} \in \que$. Then Theorem~4.3 of~\cite{mf_lf} implies that the deep hole $\btau_{k+1}$ has finite order $\ell$ in the quotient group $\real^2/\Lambda_{\tau_k}$, meaning that $\ell \btau_{k+1} \in \Lambda_{\tau_k}$. Since $\ell \be_1 \in \Lambda_{\tau_k}$, we have that
$$\ell \Lambda_{\tau_{k+1}} = \ell \spn_{\zed} \{ \be_1,\btau_{k+1} \} \subseteq \Lambda_{\tau_k} = \spn_{\zed} \{ \be_1, \btau_k \},$$
thus $\Lambda_{\tau_k}$ contains a similar copy of $\Lambda_{\tau_{k+1}}$ as a sublattice. This implies that the corresponding elliptic curves $E_{\tau_k}$ and $E_{\tau_{k+1}}$ are isogenous. Since isogeny is an equivalence relation, we have that the entire sequence of elliptic curves $\{ E_{\tau_k} \}_{k=0}^n$ is isogenous.
\smallskip

Next, we prove a bound on the smallest degree of an isogeny. Recall that $d_k = \min \{ d \in \zed_{>0} : da_k, d^2b^2_k \in \zed\}$, and notice that the fundamental deep hole of the lattice $d_k \Lambda_{\tau_k} = \spn_{\zed} \{ d_k \be_1, d_k \btau_k \}$ is $d_k \btau_{k+1}$. Successive minima of this new lattice are $d_k \lambda_{k1}$, $d_k \lambda_{k2}$; furthermore,
$$(d_k \lambda_{k1})^2,\ (d_k \lambda_{k2})^2,\ d_k \be_1 \cdot d_k \btau_k \in \zed.$$
Then Theorem~4.3 of~\cite{mf_lf} implies that the deep hole $d_k \btau_{k+1}$ has order 
$$\ell \leq 12\sqrt{3}\ (d_k \lambda_{k2})^4 = 12\sqrt{3}\ d_k^4\ |\tau_k|^4$$
in $\real^2/\Lambda_{d_k \tau_k}$, which is the same as the order of $\btau_{k+1}$ in $\real^2/\Lambda_{\tau_k}$. Then
$$[\Lambda_{\tau_k} : \ell \Lambda_{\tau_{k+1}}] = \frac{\ell^2 \det \Lambda_{\tau_{k+1}}}{\det \Lambda_{\tau_k}} = \frac{\ell^2 b_{k+1}}{b_k} \leq \frac{12\sqrt{3}\ b_{k+1}\ d_k^4\ (a_k^2+b_k^2)^2}{b_k}.$$
There exists an isogeny between $E_{\tau_k}$ and $E_{\tau_{k+1}}$ whose degree is $\delta_k = [\Lambda_{\tau_k} : \ell \Lambda_{\tau_{k+1}}]$.
\endproof
\bigskip

\section{Counting similarity classes with a prescribed deep hole}
\label{count}

In this section we will prove Theorem~\ref{ht_count}, starting with some notation. Throughout this section, let $K$ be a number field of degree $n$ over $\que$ with at least one real embedding, let $\Delta_K$ be the discriminant of $K$, and let $\O_K$ be its ring of integers. Write $\sigma_1,\dots,\sigma_{r_1}$ for the real embeddings of $K$. Specifically, we fix the real embedding $\sigma_1$ and identify every element $c \in K$ with the corresponding real number $\sigma_1(c)$. Let also $\tau_1,\bar{\tau}_1,\dots,\tau_{r_2},\bar{\tau}_{r_2}$ be the complex embeddings of $K$, then
$$n = r_1+2r_2.$$ 
The space $K_{\real} := K \otimes_{\que} \real \cong \real^{r_1} \times \cee^{r_2}$ can be viewed as a subspace of
$$\left\{ (\bx,\bwy) \in \real^{r_1} \times \cee^{2r_2} : y_{r_2+j} = \bar{y}_j\ \forall\ 1 \leq j \leq r_2 \right\} \subset \cee^n.$$
Then for any $m \geq 1$, $K_{\real}^m$ is a $mn$-dimensional Euclidean space under the bilinear form induced by the trace form on $K$:
$$\left< \alpha,\beta \right> := \Tr_K(\alpha \bar{\beta}) \in \que,$$
for every $\alpha,\beta \in K$, where $\Tr_K$ is the absolute number field trace on $K$. We also define the sup-norm on $K_{\real}^m$ by
$$|\bx| := \max \{ |x_1|,\dots,|x_{mn}| \},$$
for any $\bx \in K_{\real}^m$, where $|x_j|$ stands for the usual absolute value of $x_j$ on $\cee$. Let 
$$\Sigma_K = (\sigma_1,\dots,\sigma_{r_1},\tau_1,\bar{\tau}_1,\dots,\tau_{r_2},\bar{\tau}_{r_2}) : K \hookrightarrow K_{\real}$$
be the Minkowski embedding. We also write $\Sigma_K$ for the induced map $K^m \hookrightarrow K_{\real}^m$. Then $\Lambda_K^m := \Sigma_K(\O_K^m)$ is a lattice of full rank in $K_{\real}^m$. We define the determinant of a full-rank lattice to be the absolute value of the determinant of any basis matrix for the lattice, then
\begin{equation}
\label{ideal_det}
\det(\Lambda_K^m) = |\Delta_K|^{m/2},
\end{equation}
as follows, for instance, from Corollary 2.4 of \cite{bayer} (see also Theorem~2 of~\cite{thunder}).
\smallskip

Let $M(K)$ be the set of places of $K$ and let $n_v = [K_v : \que_v]$ be the local degree of $K$ at each place $v \in M(K)$. Select absolute values so that $|\ |_v$ extends the usual absolute value on $\que$ if $v$ is archimedean and $|\ |_v$ extends the usual $p$-adic absolute value on $\que$ if $v$ is non-archimedean. With this choice, the product formula reads
$$\prod_{v \in M(K)} |c|_v^{n_v} = 1,$$
for each $0 \neq c \in K$. For each archimedean place $v$, $|c|_v = |\sigma_j(c)|$ or $|\tau_j(c)|$ for an appropriate $1 \leq j \leq r_1$ or $1 \leq j \leq r_2$, depending on whether $v$ is real (i.e., $n_v=1$) or complex (i.e, $n_v=2$); here $|\ |$ stands for the usual absolute value on $\cee$. For each $v \in M(K)$, define the local sup-norm of a vector $\bx \in K^m$ as
$$|\bx|_v = \max \{ |x_1|_v,\dots,|x_m|_v \}.$$

For a point $\bx \in K^m$, define its {\it denominator} to be
\begin{equation}
\label{den}
d(\bx) = \min \{ c \in \que_{> 0} : c \bx \in \O_K^m \},
\end{equation}
and let the {\it (rationally) primitive point} corresponding to $\bx$ be $\bx_p = d(\bx) \bx$. 

\begin{lem} \label{prim1} The denominator $d(\bx)$ (and hence the primitive point $\bx_p$) is well-defined.
\end{lem}

\proof
We need to show that the minimum in the definition~\eqref{den} exists. Let 
$$\bx = \left( x_1,\dots,x_m \right) \in K^m.$$
There exists a rational number $c$ such that $c\bx \in \O_K^m$, and write $c = \frac{c_1}{c_2} \in \que$ with $\gcd(c_1,c_2)=1$. Then
$$\frac{c_1}{c_2} x_j \in \O_K,\ \forall\ 1 \leq j \leq m,$$
and hence the number field norm
$$\NN_K \left( \frac{c_1}{c_2} x_j  \right) = \frac{ c_1^n}{c_2^n} \NN_K(x_j) \in \zed,$$
while $\NN_K(x_j) = x_{j1}/x_{j2} \in \que$ with $\gcd(x_{j1},x_{j2}) = 1$. This means that $c_2 \mid x_{j1}$ for all $1 \leq j \leq m$, which means that the set of possible denominators for $c$ is finite, and hence the minimum of $c$ exists.
\endproof

Notice that for any $\bx \in K^m$ and non-archimedean $v \in M(K)$, $|\bx_p|_v \leq 1$, since $\bx_p \in \O_K^m$. We define the {\it primitive height} of $\bx \in K^m$ to be
$$\H^p(\bx) := \max_{v \mid \infty} |\bx_p|_v.$$
While it is clear that when $\bx_p = \bwy_p$, then $\bx$ and $\bwy$ represent the same projective point over $K$, it may happen that $\bx_p \neq \bwy_p$ while $\bx = \alpha \bwy$ for some $\alpha \in K$ (e.g., $\bx = -\bwy$).
\smallskip

\begin{figure}[t]
\centering
\includegraphics[scale=0.4]{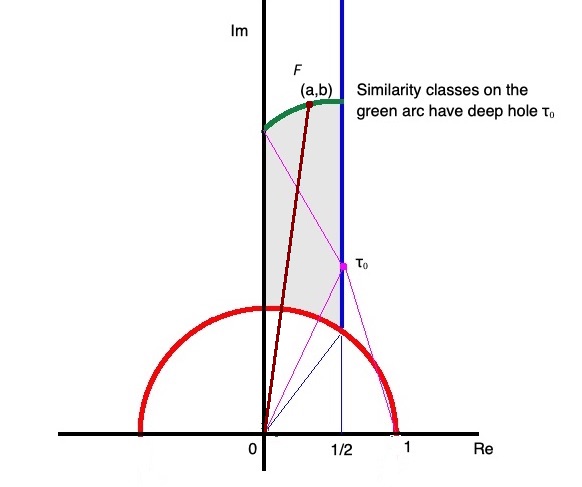}
\caption{Similarity classes with a prescribed deep hole. Pink lines are radii of the circle centered at $\tau_0$. The brown line $y = \frac{b}{a} x$ intersects the green arc at a point $\tau = a+bi$ defined over $K$.}\label{fig:dh_sim}
\end{figure}

Now, let $\tau_0$ be the fundamental deep hole lying over $K$ for some lattice $\Lambda_{\tau}$, $\tau = a+bi \in \F$. Then $\Re(\tau_0) = 1/2$ and \eqref{dh_seq1} implies that
$$t := \Im(\tau_0) = \frac{a^2+b^2-a}{2b} \in K,\ b > t.$$ 
The set of $\tau \in \F$ so that $\tau_0$ is the deep hole for $\Lambda_{\tau}$ can then be described as
\begin{equation}
\label{green_arc}
S_{\tau_0} := \left\{ a+bi :  0 \leq a \leq \frac{1}{2},\ b > t,\ \left( a- \frac{1}{2} \right)^2 + (b-t)^2 = \frac{1}{4} + t^2 \right\},
\end{equation}
which is precisely the green arc in Figure~\ref{fig:dh_sim}. Notice that each line $y=sx$ through the origin with slope 
\begin{equation}
\label{line_slope}
s \geq 2t+\sqrt{4t^2+1}
\end{equation}
intersects $S_{\tau_0}$ in precisely one point. Conversely, every point on $S_{\tau_0}$ is the intersection point with a unique such line. Hence,  $S_{\tau_0}$ is in bijective correspondence with the set of slopes $[2t+\sqrt{4t^2+1}, \infty)$. Further, notice that $\bo$ is a point on the circle of which $S_{\tau_0}$ is an arc. Let $K$ be a number field, and notice the line $y=sx$ through $\bo$ intersects $S_{\tau_0}$ in a point $(a,b) \in K^2$ if and only if $s \in K$. This means that the set $S_{K,\tau_0}$ as defined in~\eqref{SK} is given by $\left\{ a+bi \in S_{\tau_0} : a,b \in K \right\}$, and hence is in bijective correspondence with the set of slopes
$$\M_{K,\tau_0} := \left\{ s \in K : s \geq 2t+\sqrt{4t^2+1} \right\}.$$
For each $s \in \M_{K,\tau_0}$ there exists a unique $a+bi \in S_{\tau_0}$ such that $s=b/a$, and for each such pair $(a,b)$ there is a primitive point $(\alpha,\beta) \in \O_K^2$ such that $b/a = \beta/\alpha$, i.e. $(\alpha,\beta)$ is the unique primitive point on the line $y=sx$. Define
$$P_{K,\tau_0} = \left\{ (\alpha,\beta) \in \O_K^2 : (\alpha,\beta) \text{ is primitive}, \beta/\alpha \geq 2t+\sqrt{4t^2+1} \right\},$$
then we have a bijective correspondence:
$$\tau = a+bi \in S_{K,\tau_0}\ \longleftrightarrow\ (\alpha,\beta) \in P_{K,\tau_0} \text{ such that } \beta/\alpha = b/a.$$
Let us write $\H^p(\tau)$ for the primitive height $\H^p(a,b)$, then for any $T \geq 1$,
\begin{equation}
\label{S_to_P}
\left| \left\{ \tau \in S_{K,\tau_0} : \H^p(\tau) \leq T \right\} \right| = \left| \left\{ (\alpha,\beta) \in P_{K,\tau_0} : |\Sigma_K(\alpha,\beta)| \leq T \right\} \right|,
\end{equation}
which is finite.

\begin{lem} \label{count_prim} Let $T \geq 1$ and let 
$$P_{K,\tau_0}(T) = \left\{ (\alpha,\beta) \in P_{K,\tau_0} : |\Sigma_K(\alpha,\beta)| \leq T \right\}.$$
Then, as $T \to \infty$,
$$|P_{K,\tau_0}(T)| \leq \left( \frac{4^{r_1} \pi^{2r_2}}{8 \zeta(2n) \left( 2t+\sqrt{4t^2+1} \right) |\Delta_K|} \right) T^{2n} +  O(T^{2n-1}),$$
where $\zeta$ stands for the Riemann zeta-function.
\end{lem}

\proof
Consider the lattice $\Lambda_K^2 \subset K_{\real}^2$ and write $(\Lambda_K^2)_{\pr}$ for the set of primitive points in $\Lambda_K^2$, i.e. points that are not integer multiples of any other points in $\Lambda_K^2$. Notice that if $(\alpha,\beta) \in P_{K,\tau_0}$, then $\Sigma_K(\alpha,\beta) \in (\Lambda_K^2)_{\pr}$. Let us write points in $K_{\real}^2$ as $(\bx_1,\bwy_1,\bx_2,\bwy_2)$, where $(\bx_1,\bwy_1), (\bx_2,\bwy_2) \in K_{\real}$. Write $x_{11}$ and $x_{21}$ for the first coordinates of $\bx_1$ and $\bx_2$, respectively, and define
$$C_{K,t}(T) = \left\{ \bz := (\bx_1,\bwy_1,\bx_2,\bwy_2) \in K_{\real}^2 : |\bz| \leq T, x_{21} \geq \left( 2t+\sqrt{4t^2+1} \right) x_{11} \geq 0 \right\},$$
then $\Sigma_K(P_{K,\tau_0}(T)) \subseteq C_{K,t}(T) \cap (\Lambda_K^2)_{\pr}$. Write $q := \left( 2t+\sqrt{4t^2+1} \right) $ and observe that $C_{K,t}(T)$ is a direct product of $2r_2$ circular discs of radius $T$, given by
$$|y_{jk}| \leq T,\ \forall\ j=1,2 \text{ and } k=1,\dots,r_2,$$
$2r_1-2$ intervals given by
$$|x_{jk}| \leq T,\ \forall\ j=1,2 \text{ and } k=2,\dots,r_1,$$
and the set $\{(x_{11},x_{21}) : 0 \leq x_{11} \leq x_{21}/q \}$. Then
$$\Vol_{2n} C_{K,t}(T) = \left( \pi T^2 \right)^{2r_2} (2T)^{2r_1-2} \left( \frac{T^2}{2q} \right) = \left( \frac{4^{r_1} \pi^{2r_2}}{8q} \right) T^{2n}.$$
Since $C_{K,t}(T) = T C_{K,t}(1)$, Theorem~2 on p.128 of \cite{lang} combined with~\eqref{ideal_det} asserts that
\begin{eqnarray}
\label{lattice_points}
\left| C_{K,t}(T) \cap \Lambda_K^2 \right| & \leq & \left( \frac{\Vol_{2n} C_{K,t}(1)}{\det \Lambda_K^2} \right) T^{2n} + O(T^{2n-1}) \nonumber \\
& = & \left( \frac{4^{r_1} \pi^{2r_2}}{8 \left( 2t+\sqrt{4t^2+1} \right) |\Delta_K|} \right) T^{2n} +  O(T^{2n-1}),
\end{eqnarray}
as $T \to \infty$.
\smallskip

Now, let $\bx_1,\dots,\bx_{2n}$ be a basis for $\Lambda_K^2$, in particular $\Lambda_K^2 = \spn_{\zed} \{\bx_1,\dots,\bx_{2n}\}$. Let $\varphi : \Lambda_K^2 \rightarrow \zed^{2n}$ be given by
$$\bwy = \sum_{j=1}^{2n} a_j \bx_j \mapsto \ba = (a_1,\dots,a_{2n})^{\top},$$
then $\bwy \in \Lambda_K^2$ is primitive if and only if $\varphi(\bwy) \in \zed^{2n}$ is primitive. Further, $\varphi$ extends to a linear map $\varphi : K_{\real}^2 \rightarrow \real^{2n}$ and the set $\varphi(C_{K,t}(T))$ is a convex polyhedron with $\bo$ on its boundary. In particular, if a lattice point $\ba \in \zed^{2n}$ is contained in $\varphi(C_{K,t}(T))$ then there is a corresponding primitive point $\ba' \in \zed^{2n}$ contained in $\varphi(C_{K,t}(T))$ such that $\ba = c\ba'$ for some $c \in \zed$. Hence, if $T' = \max \{ |\ba| : \ba \in \varphi(C_{K,t}(T)) \}$, then for a randomly chosen point in $\varphi(C_{K,t}(T))$ the probability that it is primitive is the same as for a randomly chosen point in $\{ \ba \in \zed^{2n} : |\ba| \leq T' \}$. By Theorem~1 of~\cite{nymann}, this probability tends to $1/\zeta(2n)$ as $T' \to \infty$. Since $T' \to \infty$ as $T \to \infty$, this implies that the proportion of primitive points among all lattice points in $C_{K,t}(T) \cap \Lambda_K^2$ tends to $1/\zeta(2n)$ as $T \to \infty$. Combining this observation with~\eqref{lattice_points} completes the proof of the lemma.
\endproof
\medskip

\proof[Proof of Theorem~\ref{ht_count}] The theorem follows by combining~\eqref{S_to_P} with Lemma~\ref{count_prim}.
\endproof
\bigskip

\noindent
{\bf Conflict of interest statement:} There are no conflicts of interest to be reported.
\smallskip

\noindent
{\bf Data availability statement:} Data sharing not applicable to this article as no datasets were generated or analyzed during the current study.
\medskip

\noindent
{\bf Acknowledgement:} We are very grateful to the anonymous reviewers for their helpful comments which improved the quality of this paper.
\bigskip

\bibliographystyle{plain}  

\end{document}